\documentclass[11pt, leqno]{amsart}
\usepackage{amsthm, amsmath}
\usepackage{amssymb} 
\usepackage{amscd}
\usepackage[curve,matrix,arrow,frame,tips]{xy}
\usepackage{verbatim}
\newtheorem{thm}{Theorem}[section]
\newtheorem{prop}[thm]{Proposition}
\newtheorem{lemma}[thm]{Lemma}

\newtheorem{Remarknumb}[thm]{Remark}

\newtheorem{Remark}[thm]{Remark}

\newtheorem{cor}[thm]{Corollary}

\newcounter{ex}[section]

\newcommand{\E}{{\mathcal E}}

\newcommand{\C}{{\bf C}}

\newcommand{\Q}{{\bf Q}}

\newcommand{\Gg}{{\mathcal G}}
\newcommand{\Hh}{{\mathcal H}}

\newcommand{\Z}{{\bf Z}}
\newcommand{\F}{{\mathcal F}}
\newcommand{\ti}{\tilde}
\newcommand{\Spec}{{\rm Spec }\, }
 \renewcommand{\O}{{\mathcal O}}
\renewcommand{\SS}{{\mathcal S}}

\newcommand{\Td}{{\rm Td}}

\renewcommand{\L}{{\mathcal L}}
\newcommand{\K}{{\mathcal K}}

\newcommand{\Kr}{{\rm K}}

\newcommand{\Pic}{{\rm Pic}}

\newcommand{\td}{{\mathfrak {Td}}}
\newcommand{\ch}{{\mathfrak {s}}}

\newcommand{\ctd}{{\mathfrak {CT}}}
\newcommand{\CH}{{\rm CH}}

\def\thfill{\null\nobreak\hfill}

\def\endproof{\thfill\vbox{\hrule
  \hbox{\vrule\hbox to 5pt{\vbox to 5pt{\vfil}\hfil}\vrule}\hrule}}

\begin{document}

\title[Grothendieck-Riemann-Roch]{ Integral Grothendieck-Riemann-Roch theorem }
\author[G. Pappas]{G. Pappas}
%\thanks{*Partially supported by  NSF Grant DMS05-01409. }
%\thanks{*email: pappas@math.msu.edu}
\address{Dept. of
Mathematics\\
Michigan State
University\\
E. Lansing\\
MI 48824-1027\\
USA\\
 \ \ \  {\rm \phantom{aaaaaaaaaaaaaaaaaaaaaaaaaaaaaaaa}}   {\sl \ \ \ \phantom{aa} email:\  pappas@math.msu.edu}}
\date{\today}

\maketitle
\section{Introduction}

Let   $f: X\to S$ 
be a projective morphism between two smooth quasi-projective algebraic varieties defined over the field $k$.
If $\F$ is a vector bundle on $X$, then the hypercohomology  ${\rm R}f_*(\F)$ 
of $\F$ is 
represented by a finite complex $({\mathcal E}^i)$ of vector bundles on $S$
and we can unambiguously define  the class
$f_* [\F]:=\sum_i(-1)^i[{\mathcal E}^i]$ of $ {\rm R}f_*(\F) $ in the Grothendieck group ${\rm K}_0(S)$. The Grothendieck-Riemann-Roch 
theorem ([BS]) is the identity
\begin{equation}\label{grr}
{\rm ch}(f_*[\F] )=f_*({\rm ch}( \F) {\rm Td}(T_X)){\rm Td}(T_S)^{-1}
\end{equation}
in the Chow ring with rational coefficients $\CH^*(S)_\Q=\oplus_{n}\CH^n(S)_\Q $.
Here $\rm ch $ is the Chern character and ${\rm Td}(T_X)$, $\Td(T_S)$
stand for the Todd power series evaluated at the Chern classes 
of the tangent bundle of $X$,
respectively $S$. Since both sides of (\ref{grr})
take values in $\CH^*(S)_\Q:=\CH^*(S)\otimes\Q$,  
only information  modulo torsion about the Chern classes of $f_*[\F]$
can be obtained from this identity. 

The goal of our paper is to improve on this as follows: 
Set
\begin{equation}\label{denom}
T_m=\prod_{p }p^{\left[\frac{m}{p-1}\right]}\,
\end{equation}
where the bracket  denotes the integral part and the product is over prime numbers.
The integer $T_m$ is the denominator of the degree $m$ part of the Todd  power series; notice that $m!$ divides $T_m$ which divides $T_{m'}$ for $m'\geq m$. Write
$$
{\rm ch}=\sum_{m\geq 0} \frac{\ch_m}{m!}\ , \ \ \qquad {\rm Td}=\sum_{m\geq 0}\frac{\td_m}{T_m}\ .
$$
The numerators   ${\ch}_m$ and ${\td}_m$ of the degree $m$ parts   
 are polynomials  
with integral coefficients in the Chern classes.
%$c_i$ (with ${\rm deg}(c_i)=i$).
%(Lemma \ref{factorials}).
We show that,  when  $k$ has  characteristic zero and 
the relative dimension $d=\dim(X)-\dim(S)$ is non-negative,
the   Grothendieck-Riemann-Roch (GRR) formula actually applies to calculate 
$$
 \frac{T_{d+n}}{n!}\cdot \ch_n(f_* [\F] ) 
$$
in $\CH^n(S)$ (and not just in $\CH^n(S)$ modulo torsion). The point here is that multiplication by 
$T_{d+n}$  clears all the denominators in the codimension $n$  component
of the formula (\ref{grr}); we show that the resulting identity is indeed true
in $\CH^n(S)$. For example, when $n=1$, $\ch_1=c_1$ and $\ch_1(f_* [\F] )$ is the class of the determinant
of cohomology $\det {\rm R}f_*(\F)$ in $\Pic(S)$. In this case,
we obtain
\begin{multline}\label{iGRR}
\ \ T_{d+1}\cdot [\det {\rm R}f_*(\F)]=\\ =-r(f_*[\F])\cdot  \frac{T_{d+1}}{2} \cdot c_1(T_S)+\sum_{m=0}^{d+1} \frac{T_{d+1}}{m!\cdot T_{d+1-m}}\cdot 
f_*\left[{\ch}_m(\F)\cdot  {\td}_{d+1-m}(T_X)\right] \ \
\end{multline}
%\begin{equation*}\label{iGRR}
%T_{d+1}\cdot [\det {\rm R}f_*(\F)] = \frac{ r(\F)\cdot T_{d+1}}{2} \cdot c_1(T_S) +\sum_{m=0}^{d+1}  \frac{T_{d+1}}{m!\cdot T_{d+1-m}}\cdot 
%f_*\left[{\ch}_m(\F)\cdot  {\td}_{d+1-m}(T_X)\right] 
%\end{equation*}
in ${\rm Pic}(S)=\CH^1(S)$. Here $r(f_*[\F])=\ch_0(f_*[\F])=\sum_{i} (-1)^i {\rm rank}_{\O_S} ({\rm R}^if_*(\F))$
is the (virtual) rank   of ${\rm R}f_*(\F)$ on $S$.
Note  that the ratios $T_{d+1}/(m!\cdot T_{d+1-m})$ and $T_{d+1}/2$ are  integers (Lemma \ref{factorials}). 

A few isolated cases of this result were already known:
When $f$ is a relative curve, i.e when $d=1$, 
Mumford has shown, using the moduli space of curves, that the GRR formula 
applies to calculate $12\cdot [\det {\rm R}f_*(\F)]$ in ${\rm Pic}(S)$. Since $T_{2}=12$, this also follows from (\ref{iGRR}) above,
which then generalizes Mumford's result to higher dimensions
(but only in characteristic $0$). In fact, by applying the integral GRR
formula for $d=1$ and $n\geq 2$, we also obtain new
integral relations among the pull-backs of tautological classes on $S$.
These relations were known before only up to torsion and refine
corresponding
(known or conjectural) equations in the integral cohomology of 
the mapping class group (see \S \ref{2d3}).
%For example, if $f: X\to S$ is a relative surface,
%we obtain
%\begin{equation*}
%24\cdot [\det {\rm R}f_*(\O_X)]=f_*\left[c_1(T_X)\cdot %c_2(T_X)\right]-12\chi(\O_X)\cdot c_1(T_S)\  
%\end{equation*}
%in $\CH^1(S)$. 
When the morphism $f$ is a finite \'etale cover, we obtain the Riemann-Roch theorem 
for covering maps of Fulton-MacPherson.
If $f: A\to S$ 
is an abelian scheme of relative dimension $g$, then our integral GRR formula implies
that the top Chern class of the Hodge bundle over $S$ in $\CH^g(S)$
is  annihilated by the integer $T_{2g}$. This also follows from a (stronger) 
result of 
Ekedahl-van der Geer. %  ([EvdG, Theorem 3.5]). 
See \S \ref{corol} for some more corollaries
for families of surfaces.

Our approach was inspired by the classical work of Washnitzer [W] and Fulton [F2]
on characterizing the arithmetic genus and by certain  
constructions in the theory of algebraic cobordism of 
Levine and Morel [LM]. The crucial ingredients are  Hironaka's resolution of singularities
and the weak factorization theorem  for birational maps of [AKMW]
(this is the only ingredient of our proof that has not been available
for a long time);
the use of these restricts the result to characteristic $0$. 

Here is an outline of the proof: A refinement of the classical arguments 
shows that the integral GRR identity holds 
for  projective bundles, for closed immersions and for blow ups
along smooth centers. However,   
contrary to what happens in Grothendieck's approach, 
the general result  does not follow easily from these special cases:
Indeed, the use of a projective bundle of dimension higher
than that of the variety introduces additional denominators.
To show the integral GRR identity in general, we first 
assume that $\F$ is the structure sheaf $\O_X$: We then prove that if $X$, $X'$ are smooth linearly equivalent divisors in $W\to S$ with $W$ smooth, then the integral GRR formula  holds for $X$ if and only if it holds for $X'$.  In fact, we can extend both sides of the integral GRR formula to general Weil divisors on $W$ and show that  each side respects linear equivalence. We then observe that if the result holds for $X\to S$ then it holds for a projective bundle ${\bf P}(\E)\to X\to S$. We also prove, using the factorization theorem and the result for blow ups,
that the integral GRR formula holds for $X\to S$ if and only if it holds for  any $X'\to S$
which is birationally equivalent to $X$ over $S$. Now to actually prove the formula
for $\F=\O_X$ we argue by double induction, first 
on $n$ and then on the relative dimension $d$. The result for $n=0$
is given by the Hirzebruch-Riemann-Roch theorem since $\CH^0(S)=\Z$ is  torsion-free.
%Assuming the result for all values less than $n$  we show the result for $n$. 
We then observe that when $f: X\to S$ is not dominant the result  follows
from the induction hypothesis on $n$ using resolution of singularities, factorization, 
and integral Riemann-Roch for closed immersions.
To show the result for $f$ dominant, we apply induction on $d$. 
Since $X$ is birational to the desingularization $Y'$ of a hypersurface 
$\hat Y$ in ${\bf P}^{d+1}\times_kS$, by the above, it is enough to deal with 
such a desingularization. The hypersurface $\hat Y$ is linearly equivalent
to a sum ${\rm pr}_1^{-1}(H)+{\rm pr}_2^{-1}(T_1)-{\rm pr}_2^{-1}(T_2)$, where $H$ is a smooth
hypersurface in ${\bf P}^{d+1}_k$ and $T_1$, $T_2$ are smooth divisors
on $S$; we eventually reduce to checking the identity in the simple cases
that
$X={\rm pr}_1^{-1}(H)$, or that $X={\rm pr}_1^{-1}(T_i)$, $i=1,2$.  The induction hypothesis on $d$ implies that
 the exceptional locus of the desingularization
we employ do not contribute an error to the formula: Indeed,
the components of the exceptional locus are birational to projective bundles over varieties
of smaller dimension.  The argument shows more or less simultaneously
that the result is  true when 
$\F$ is a line bundle on $X$. The case that $\F$ is a general vector bundle 
follows by using a result of Kleiman which allows us to split $\F$ after a blow-up. 

In fact, it turns out that an important part of the proof  
%of the integral GRR formula 
can also be presented as an application 
of   the ``generalized degree formula'' in the theory of algebraic cobordism 
(see Remark \ref{GenDeg}). This gives a  
somewhat different route toward the main result. 
We   chose the  
direct and classical argument above to make the paper more self-contained 
and accessible. However, we feel that this observation establishes an interesting 
connection which could be important in future developments.

Finally, let us mention that we expect that this  modification of Grothendieck's argument can be applied to the proof
of other Riemann-Roch type theorems and should also produce versions that capture torsion
information. For example, one could attempt to revisit the ``functorial" Riemann-Roch 
of Deligne ([D]) and Franke (unpublished) or Gillet's Riemann-Roch theorem for higher algebraic $\rm K$-theory  ([G]).
%\begin{comment}

%It could also  be crucial in attempting to extend this type of result to other set-ups
%(``functorial" Riemann-Roch, Riemann-Roch for higher algebraic $\rm K$-theory, etc.). 
%\end{comment}

\smallskip
\medskip

\noindent{\bf Acknowledgments:} The author would like to thank P. Deligne for a useful discussion, B. Totaro for his comments
and the Institute for Advanced Study for its hospitality during the year
2004-2005.
\bigskip
%\vfill\eject
 
\section{Preliminaries} \label{factorS}

\setcounter{equation}{0}
Throughout the paper $k$ is a field of characteristic $0$; 
all algebraic varieties and morphisms are over the field $k$.

\subsection{} We start with the following lemma which will be used repeatedly.

\begin{lemma} \label{factorials}
Let $m$ be a positive integer. If  $m_1+m_2+\cdots +m_r+m_{r+1}+\cdots +m_{r+s}\leq m$
with $m_i$ positive integers, then the product
\begin{equation*}
(m_1+1)!\cdots (m_r+1)!\cdot T_{m_{r+1}}\cdots \cdot T_{m_{r+s}}
\end{equation*}
divides $T_m$.
\end{lemma}

\begin{Proof}
Recall that if $p$ is a prime number and $n$ an integer  with $p^k\leq n< p^{k+1}$, then the largest power of $p$ that divides $n!$
is 
$$
 \left[\frac{n}{p}\right]+ \left[\frac{n}{p^2}\right]+\cdots +\left[\frac{n}{p^k}\right]
 \leq \left[\frac{n(p^k-1)}{p^k(p-1)}\right]\leq  {\left[\frac{n-1}{p-1}\right]}\, .
$$
The lemma now follows from this and (\ref{denom}).\endproof
\end{Proof}

\begin{comment}
 Using (\ref{denom}) and this fact, we see that the product $(m_1+1)!\cdots (m_r+1)!\cdot T_{m_{r+1}}\cdots T_{m_{r+s}}$
divides
\begin{equation*}
\prod_pp^{\left[\frac{m_1}{p-1}\right]+\cdots +\left[\frac{m_{r+s}}{p-1}\right]}
\end{equation*}
The lemma now follows from (\ref{denom}).
\end{comment}

\subsection{} Consider the  Todd power series 
\begin{equation*}
\Td=\prod_{j=1}^\infty \frac{x_j}{1-e^{-x_j}}=1+\frac{1}{2}c_1+\frac{1}{12}(c^2_1+c_2)+\frac{1}{24}c_1c_2+\cdots %-\frac{1}{720}(c^4_1-4c_1^2c_2-3c^2_2-c_1c_3+c_4)+\cdots \ %\in \Q[[c_1,c_2, c_3, \ldots ]],
\end{equation*}
 viewed
as  a formal power series with rational coefficients in the variables $c_i$ (the elementary symmetric functions of $x_j$)
with ${\rm deg}(c_i)=i$. For any $m$,
we will consider  the degree $m$ part of $\Td$ which we will denote by $\Td_m$. 
(In general, we will denote by $P_m$ the homogeneous degree $m$ part of
$P$.) By [Hi, Lemma 1.7.3], the polynomial  $$\td_{ m}=T_m\cdot \Td_{ m}$$
has integral coefficients and is the numerator of the degree $m$ part of $\Td$.
 
We can also consider the Chern power series 
$$
{\rm ch}=r+\sum_{j=1}^\infty(e^{x_j'}-1)=r+c'_1+\frac{1}{2}(c'^2_1-2c'_2)+\frac{1}{6}(c'^3_1-3c'_1c'_2+3c'_3)+\cdots
%\frac{1}{24}(c^4_1-4c^2_1c_2+4c_1c_3+2c^2_2-4c_4)+\cdots \ %\in \Q[[r,c'_1, c_2',\, \ldots ]]
$$
as a formal power series with rational coefficients in the variables $r$ (rank), $c'_i$ (the elementary symmetric functions of $x'_j$)
with $\deg(r)=0$, ${\rm deg}(c'_i)=i$.  
We set $\ch_m=m!\cdot {\rm ch}_m$ for the numerator 
of the degree $m$ part of ${\rm ch}$.
Also set 
$$
\ctd_m=T_{m}\cdot ({\rm ch}\cdot \Td)_{m}=\sum_{j=0}^m\frac{T_m}{j!\cdot T_{m-j}}\cdot (\ch_j\cdot \td_{m-j})\ .
$$
By Lemma \ref{factorials}, $\ctd_{m}$ is a  homogeneous  polynomial  in $\Z[c_1,c_2,\ldots , c_{m}, r, c'_1,\ldots, c'_m]$.

\subsection{} Let $Y$ be a variety over $k$. We will denote by ${\rm K}_0(Y)$ the Grothendieck ring
of locally free coherent $\O_Y$-sheaves on $Y$ and by ${\rm G}_0(Y)$ the Grothendieck 
group of coherent $\O_Y$-sheaves on $Y$. Suppose that $Y$ is smooth and quasi-projective.
Then
the natural map ${\rm K}_0(Y)\to {\rm G}_0(Y)$ is an isomorphism;
we will   identify these two groups without further notice. 
Denote   by $\CH^i(Y)$ the Chow
group of algebraic cycles of codimension $i$ on $Y$
modulo rational equivalence. 
There are well-defined intersection pairings
$\CH^i(Y)\otimes \CH^j(Y)\to \CH^{i+j}(Y)$
which turn $\CH^*(Y)=\oplus_{i=0}^{\dim(Y)}\CH^i(Y)$ into 
a graded commutative ring.  %Set $\CH^{\leq m}(Y)=\oplus_{0\leq i\leq m}\CH^i(Y)$.
If $\F$ is a locally free coherent $\O_Y$-sheaf on $Y$ we have the Chern classes $c_i(\F)\in \CH^i(Y)$,
$1\leq i\leq \dim(Y)$. We will denote by
$T_Y:=(\Omega^1_{Y/k})^{\vee}$ the tangent sheaf of $Y$.
For $m\geq 0$, we now set
\begin{eqnarray*}
&&\ \ \ \ctd_m(\F,Y):=\ctd_{ m} (c_1(T_Y),   \ldots, c_{m}(T_Y),   r(\F), c_1(\F),\ldots, c_{m}(\F)),\\
 &&\ \ \ \ \ \ \  \ \ \td_{m}(\F) :=\td_{ m} (c_1(\F),\ldots, c_{m}(\F)),
\end{eqnarray*}
 in $\CH^m(Y)$. (We evaluate $\ctd_{m}$  by setting 
$c_i=c_i(T_Y)$, $r={\rm rank}(\F)$, $c'_i=c_i(\F)$, similarly for $\td_{ m}(\F )$.)
%There are also non-homogeneous versions $\ct_{\leq m}(\F; X/S)$,  $\td_{\leq m}(\F; X/S)$, 
%$\td^0_{\leq m}(\F; X/S)$. 
Often, we will simply write $\td_{m}(Y)$ 
instead of $\td_{m}(T_Y)$. 
It follows from the Whitney sum formula
that the functions $\ch_m(-)$, $\td_m(-)$ and
$\ctd_m(-,Y)$ extend to give well-defined maps ${\rm K}_0(Y)\to \CH^{ m}(Y)$. The maps $\ch_m(-)$ and $\ctd_m(-,Y)$ 
are additive. The multiplicativity of the Chern, resp. Todd, power series
implies
\begin{equation}\label{mulChern}
\ \ \ch_m(a\cdot b)=\sum_{i=0}^m \frac{m!}{i!\cdot (m-i)!}  \cdot \ch_i(a)\cdot \ch_{m-i}(b)\,,
\end{equation}
\begin{equation}\label{mulTodd}
\td_m(a+b)=\sum_{i=0}^m \frac{T_m}{T_i\cdot T_{m-i}} \cdot \td_i(a)\cdot \td_{m-i}(b)\, ,
\end{equation}
with $a$, $b$ in the Grothendieck ring  ${\rm K}_0(Y)$.

%If in addition, 
%$m=\dim(X)-\dim(S)+1=d+1$, we will simply write $\td (X/S)$ 
%instead of $\td_{d+1}(X/S)$.

Suppose that $f: X\to S$ is a projective morphism between
the smooth varieties $X$ and $S$. Set $d=d_f=\dim(X)-\dim(S)$. There are   well-defined push-forward homomorphisms:
\begin{equation*}
f_*: \CH^i(X)\xrightarrow{ \ }\CH^{i-d}(S),\qquad
%\end{equation*} 
%\begin{equation*}
f_*: {\rm K}_0(X)= {\rm G}_0(X)\xrightarrow{ \ }{\rm K}_0(S)= {\rm G}_0(S) ,
\end{equation*}
where for $\F$  a (locally free) coherent $\O_X$-sheaf, we set
$f_*[\F]=[{\rm R}f_*(\F)]=\sum_{i}(-1)^i[{\rm R}^if_*(\F)]$.

Our main result is:

\begin{thm}\label{main}
Suppose $k$ is a field of characteristic $0$. Let $X$ and $S$ be 
smooth quasi-projective varieties
over $k$ and let $f: X\to S$ be a projective morphism 
over $k$.  Set $d=d_f=\dim(X)-\dim(S)$ and suppose that $\F$ is a 
coherent   $\O_X$-sheaf   on $X$. 

a) Suppose $d\geq 0$. Then the identity
\begin{equation}\label{RReqA}
\frac{T_{d+n}}{T_n}\cdot  \ctd_{ n}(f_*[\F], S)=f_*(\ctd_{d+n}(\F, X))\ \ \ \ \
\end{equation}
holds in $\CH^n(S)$.

b) Suppose $d<0$. Then the identity
\begin{equation}\label{RReqB}
\ \ \ \ \ \ctd_{ n}(f_*[\F], S)=\frac{T_{n}}{T_{n+d}}\cdot f_*(\ctd_{n+d}(\F, X))
\end{equation}
holds in $\CH^n(S)$.
\end{thm}

Let $\ctd_{m}(\F, X/S)$ in $\CH^m(X)$ be the result of evaluating the polynomial
$\ctd_m$ by setting $c_i=c_i([T_X]-[f^*T_S])$, $r={\rm rank}(\F)$, $c'_i=c_i(\F)$. Part (a) implies the following.
 
\begin{cor}\label{maincor}
Suppose $d\geq 0$. Then the identity
\begin{equation}\label{RReq}
\frac{T_{d+n}}{n!}\cdot  \ch_n(f_*[\F])=f_*(\ctd_{d+n}(\F, X/S))
\end{equation}
holds in $\CH^n(S)$.
\end{cor}

\begin{Proof}
Observe that the left hand side of the identity (\ref{RReqA}) can be written
\begin{equation}
 \frac{T_{d+n}}{n!}\cdot \ch_n(f_*[\F])+\sum_{j=1}^{n}\frac{T_{d+n}}{T_{d+n-j}\cdot T_j}\cdot
 \left\{ \frac{T_{d+n-j}}{(n-j)!}\cdot \ch_{n-j}(f_*[\F])\right\}\cdot \td_{j}(T_S)\, .
\end{equation}
The statement follows from this observation, the projection formula and (\ref{RReqA}), by induction on $n$.\endproof
\end{Proof}

\begin{Remark}\label{bigRem}
{\rm  a) The image of $f_*(\ctd_{d+n}(\F, X/S))$ in $\CH^{n}(S)\otimes\Q$ is
$$
T_{d+n}\cdot f_*(({\rm ch}(\F)\cdot {\rm Td}(T_X)\cdot \Td(f^*T_S)^{-1})_{d+n})
$$
and so the image of the identity (\ref{RReq}) in ${\CH}^n(S)\otimes\Q$
is the identity for  ${\rm ch}_n( f_*[\F])$ given by the  
Grothendieck-Riemann-Roch theorem.

\begin{comment}
b) The identity (\ref{iGRR}) of the introduction follows easily from  
(\ref{RReq}) for $n=1$ and the Hirzebruch-Riemann-Roch theorem applied to the generic fiber
of $f: X\to S$.

c) Given a morphism $f$ of non-negative relative dimension, 
a sheaf $\F$ and an integer $n\geq 0$, we can  see that if (\ref{RReq}) holds for 
$f$, $\F$ and all values $n'\leq n$, then (\ref{RReqA}) also holds for $f$, $\F$ and $n$.
\end{comment}

b) When the morphism $f$ is a closed immersion of codimension $r=-d$,
then (\ref{RReqB}) follows from the ``Riemann-Roch without denominators" of Jouanolou [J] (see Theorem \ref{JoRR}).
 }
 \end{Remark}

\subsection{} \label{corol}
Here we describe some corollaries of this result.

\subsubsection{} Let $f: X\to Y$ be a finite \'etale morphism
between smooth quasi-projective varieties over $k$. Then $f^*T_Y\simeq T_X$.
Therefore, (\ref{RReq}) implies
\begin{equation*}\label{FM1}
\frac{T_n}{n!}\cdot \left(\ch_n(f_*\F)-  f_*(\ch_n(\F))\right)=0\, ,
\end{equation*}
for any $\F$ on $X$. As in [FM, Remark 23.8], we see that
this immediately implies
\begin{equation}\label{FM}
L_n\cdot \left(\ch_n(f_*\F)-  f_*(\ch_n(\F))\right)=0\, ,
\end{equation}
where $L_n$ is the product of all primes that divide $T_n/n!$.
(The integer ${T_n}/{n!}$
is denoted by $N_n$ in loc. cit.)
This last identity (\ref{FM}) is the integral Riemann-Roch theorem
for covering maps of Fulton-MacPherson 
([FM, Theorem 23.3]). In the context of group representations and for characteristic classes in integral group cohomology,  Evans-Kahn [EK] have shown that, for the (topological) cover
given by $BH\to BG$ where $H$ is a subgroup of a finite group
$G$,  the integers $L_n$ are the smallest with the property 
corresponding to (\ref{FM}).
Using Totaro's construction [T], we can approximate $BH\to BG$ by a finite \'etale cover of smooth quasi-projective varieties $X\to Y$. Hence, we see that $L_n$ are the smallest integers so that (\ref{FM}) holds for all finite \'etale covers.

\subsubsection{}  Let $f: A\to S$ be an abelian scheme ([CF]) of relative dimension $g$ 
over the smooth quasi-projective variety $S$ over  
$k$. By a result of Grothendieck the morphism $f$ is projective.
Using (\ref{RReq}) we obtain
\begin{equation}\label{evdg2}
 \frac{T_{2g}}{g!}\cdot \ch_g(f_*[\O_A])= f_*(\ctd_{2g}(\O_A, A/S)) 
\end{equation}
in $\CH^g(S)$.
 The Hodge bundle
is the locally free coherent $\O_S$-sheaf $E=s^*(\Omega^1_{A/S})$
where   $s: S\to A$
is the zero section; it has rank $g$ and we have $\Omega^1_{A/S}\simeq f^*(E)$. 
We find 
$$
f_*(\ctd_{2g}(\O_A, A/S))=f_*(\td_{2g}(f^*(E^{\vee})))
=f_*f^*(\td_{2g}(E^{\vee}))=0
$$ 
in $\CH^g(S)$, %Set $\lambda_g=c_g(E)$ in $\CH^g(S)$. 
while $f_*[\O_A]=\sum_{i=0}^g(-1)^i[{\rm R}^if_*(\O_A)]=\sum_{i=0}^g(-1)^i[\wedge^i(E^\vee)]$. 
The standard identity [BS, Lemme 18] now gives
\begin{equation}\label{72ab}
 \ch_g\left(\sum_{i=0}^g(-1)^i[\wedge^i(E^\vee)]\right)
= g!\cdot c_g(E)\, .
\end{equation}
Therefore, (\ref{evdg2}) implies that $T_{2g}\cdot c_g(E)=0$   in $\CH^g(S)$.
For $g=1$, we get the classical $12\cdot c_1(E)=0$. 
Ekedahl and van der Geer show that $2(g-1)!\,D_{2g}\cdot c_g(E)=0$ ([EvdG,  Theorem 3.5])
where
$$
D_{2g}= \prod_{l\,  {\rm \, prime,}\,\, l-1|2g}l^{1+{\rm ord}_l(2g)}\ .
$$
By von Staudt's theorem, the number $D_{2g}$ is equal to   the denominator of $B_{2g}/2g$ with $B_{2g}$ the Bernoulli number.  
When $g>1$, we can see that $2(g-1)!\,D_{2g}$ divides $T_{2g}$ and so this corollary of (\ref{RReq})  follows from their result.  

\subsubsection{}\label{2d3}
Assume in addition that $f: X\to S$ is smooth and that the geometric fibers of $f$ are   irreducible curves ($d=1$). Let 
$$
\kappa_{i}=f_*(c_1(\Omega^1_{X/S})^{i+1})
$$ 
be Mumford's ``tautological" classes in $\CH^i(S)$. Denote by $\omega={\rm R}^0f_*(\Omega^1_{X/S})$ the Hodge bundle on $S$. Since ${\rm R}^0f_*(\O_X)\simeq \O_S$ and ${\rm R}^1f_*(\O_X)\simeq \omega^\vee$
(by Serre-Grothendieck duality), we have
$\ch_n(f_*[\O_X])= (-1)^{n-1}\ch_n(\omega)$. Applying (\ref{RReq}) to $\F=\O_X$ 
now gives  $12\cdot c_1(\omega)=\kappa_1$ in ${\rm Pic}(S)$ (for $n=1$) and 
\begin{equation}\label{MMclass}
\frac{T_{n+1}}{n!}\cdot \ch_n(\omega)=\begin{cases}
\ 0,&  \hbox{\rm if $n=2m$}\,, \\
\displaystyle{T_{2m}\frac{B_{2m}}{(2m)!}\cdot  \kappa_{2m-1}},&  \hbox{\rm if $n=2m-1$}\,,
\end{cases}
\end{equation}
in $\CH^{n}(S)$ for $n\geq 2$. %(Here $B_{2m}$ is the Bernoulli number.) 
The corresponding identity in the  integral  cohomology of the mapping class group of surfaces is a slightly weakened version of a conjecture of Akita [A].
A corollary of (\ref{MMclass}) is 
that we have
$
\kappa_{2m-1}=D_{2m}\cdot \alpha_m+\beta_m
$
in $\CH^{2m-1}(S)$, where $D_{2m}$ is the denominator of $B_{2m}/2m$ as above,
and $\beta_m$ is ${T_{2m}}/({D_{2m}\cdot (2m-1)!})$-torsion.
In fact, we conjecture that $\beta_m$ can be taken to be zero,
i.e that $\kappa_{2m-1}$ is actually
$D_{2m}$-divisible in $\CH^{2m-1}(S)$.
(If $k=\C$, see [GMT] for a discussion of the corresponding statements in the 
integral cohomology ${\rm H}^{4m-2}(S(\C),\Z)$.)

\subsubsection{} Let $f: X\to S$ be a   relative  surface ($d=2$).
Assume that $f$ is smooth and set $\omega_{X/S}=\det(\Omega^1_X)\otimes_{\O_X} \det(f^*\Omega_S)^{-1}$ for the relative dualizing sheaf. For $m\geq 0$,
apply (\ref{RReq}) to $n=1$ and the sheaves $\O_X$ and $\omega_{X/S}^{\otimes\, m}$.
We deduce the existence of an isomorphism of invertible sheaves on $S$ 
\begin{equation*}\label{isodual}
\left(\det {\rm R}f_*(\omega_{X/S}^{\otimes\, m})\otimes_{\O_S} \det {\rm R}f_*(\O_X )^{\otimes  (2m-1)}\right)^{\otimes 24}\simeq \langle\omega_{X/S}, \omega_{X/S}, \omega_{X/S}\rangle^{\otimes\, m(6m-4m^2-2)} 
\end{equation*}
where the bracket denotes Deligne's intersection
bundle ([D]). (By loc. cit., the class of $\langle\omega_{X/S}, \omega_{X/S}, \omega_{X/S}\rangle$ in $\Pic(S)=\CH^1(S)$ is equal to $f_*(c_1(\omega_{X/S})^3)=
-f_*(c_1([T_X]-[f^*T_S])^3)$.) It would be interesting to establish a {\sl canonical} isomorphism as above. 

Let us consider an application:
Suppose $f: X\to S$ is  a    family of   Enriques surfaces. Then $\omega_{X/S}^{\otimes\, 2}$
is trivial along the fibers of $f$. Therefore, $\K:={\rm R}^0f_*(\omega_{X/S}^{\otimes\, 2})$ 
is an invertible sheaf  on $S$ and we have  
 $\omega_{X/S}^{\otimes\, 2}\simeq f^*\K$.  We also have $ \det {\rm R}f_*(\O_X )\simeq \O_S$; hence,  the projection formula gives 
$\det {\rm R}f_*(f^*\K)\simeq \K$.  The above isomorphism for $m=2$ now
gives 
$$
\K^{\otimes\, 24}\simeq \langle \omega_{X/S}, \omega_{X/S}, \omega_{X/S}\rangle^{-\otimes\, 12} \simeq \langle f^*\K, f^*\K, \omega_{X/S}\rangle ^{-\otimes\, 3}\simeq \O_S\,.
$$
By [B], such a trivialization of a power of $\K$ can be given 
explicitly using a Borcherds product on the period domain.

\begin{comment}
\subsubsection{} Now let $f: X\to S$ be a  family of Enriques surfaces,
and consider the isomorphism above for $m=2$. Since $\omega_{X/S}^{\otimes\, 2}$
is trivial along the fibers of $f$, $\K:={\rm R}^0f_*(\omega_{X/S}^{\otimes\, 2})$ 
is an invertible sheaf  on $S$ and we have  
 $\omega_{X/S}^{\otimes\, 2}\simeq f^*\K$.  We also have $ \det {\rm R}f_*(\O_X )\simeq \O_S$; hence  the projection formula gives 
$\det {\rm R}f_*(f^*\K)\simeq \K$. Therefore we obtain 
$$
\K^{\otimes\, 24}\simeq \langle \omega_{X/S}, \omega_{X/S}, \omega_{X/S}\rangle^{-\otimes\, 12} \simeq \langle f^*\K, f^*\K, \omega_{X/S}\rangle ^{-\otimes\, 3}\simeq \O_S\,.
$$
%This last corollary should also follow from the theory of Borcherds' products
%([Bo]).
%which constructs a trivialization of $\K^{\otimes 4}$ 
\end{comment}
 
\subsection{} We will  say   that
integral Riemann-Roch   holds for $(f, n)$
when either (\ref{RReqA}) or (\ref{RReqB}) (depending 
if $d_f\geq 0$ or $d_f<0$) holds for
all $\F$ on $X$. The following observation will
be used repeatedly in our proof of Theorem \ref{main}.

\begin{prop}\label{compose}
Let  $f: X\to Y$ and $g: Y\to Z$ be projective morphisms
between smooth quasi-projective varieties over $k$.
Suppose that  
integral Riemann-Roch  holds for both $(f, n+d_g)$ and $(g, n)$. Suppose
in addition that either $d_f\geq 0$ or $d_g\leq 0$. 
Then integral Riemann-Roch holds for $(g\cdot f, n)$. 
\end{prop}

\begin{Proof}
This follows easily from the fact that the push-forward homomorphisms  
(both for Grothendieck groups and Chow groups)
satisfy
$(g\cdot f)_*=g_*\cdot f_*$. (The assumption on $d_f$, $d_g$ is needed
to guarantee that certain ratios of Todd denominators which are involved in the argument are integers.)\endproof
\end{Proof}

\begin{Remark}{\rm
Grothendieck's proof of the  Riemann-Roch theorem ([BS])
involves factoring a morphism into a composition of a closed immersion $f$
followed by a projective bundle $g$. In that case, 
$d_f< 0$ and $d_g> 0$, and so Proposition \ref{compose}
does not apply. 
}
\end{Remark}

\begin{prop}\label{multiply}
Let  $f: X\to Y$  be a  projective morphism
between smooth quasi-projective varieties over $k$.
Let $\F$ be a coherent $\O_X$-sheaf on $X$ and $\Gg$ a locally free 
coherent $\O_Y$-sheaf on $Y$. Given $n\geq 0$, suppose that the integral Riemann-Roch formula holds for $f$, $\F$, and  all $n'\leq n$. Then it also holds for
$f$, $\F\otimes_{\O_X} f^*\Gg$, and all $n'\leq n$. 
\end{prop}

\begin{Proof}
The proof follows easily  from (\ref{mulChern})  and the projection formula.\endproof
\end{Proof}

\smallskip

\section{Divisors}

\setcounter{equation}{0}

\subsection{} 
For $m\geq 1$, let us consider the polynomial
$$
Q_{m}(  c_1,   \ldots , c_{m-1},   x)=T_{m-1}\cdot ((1-e^{-x}) \cdot \Td )_{m}
$$
in the variables   $c_1,  \ldots  , c_{m-1} $, $x$, with $\deg(c_i)=i$, $\deg(x)=1$.
By Lemma \ref{factorials}, $Q_{m}$ 
has integral coefficients.

Suppose that $W$ is a smooth quasi-projective variety of dimension $\delta+1\geq 1$
over $k$.
If $[D]\in \CH^1(W)$ is the class of the Weil divisor $D=\sum_in_iD_i$ of $W$,  
we can consider
\begin{equation*}
\td_{m}(D; W ):=Q_{m}(c_1(T_W),   \ldots , c_{m-1}(T_W),    [D])
\end{equation*}
in $\CH^m(W)$. Notice that, by its definition, 
 $\td_m(D;W ) $  depends only on the linear
equivalence class $[D]$ of $D$. We also have
\begin{equation}\label{32}
\frac{T_{m}}{T_{m-1}}\cdot \td_{m}(D; W )=\ctd_{m}([\O_W]-[\O_W(-D)], W )\,, 
\end{equation}
where the right hand side is defined in Section \ref{factorS}.

\begin{prop} \label{vample1}
a) Suppose that $D$ is a smooth divisor
and denote by $i: D\hookrightarrow W$ the natural embedding. 
Then we have
\begin{equation}\label{smooth}
\td_{m}(D; W )=i_*(\td_{m-1 }(D ))
\end{equation}
in $\CH^{m}(W)$.

b) Suppose that $D=D_1+D_2$ with $D_1$, $D_2$ smooth. Suppose also
that the scheme theoretic intersection $D_1\cap D_2$
is smooth and of pure codimension $2$. 
Then
\begin{equation}
\td_{m}(D;W )  =(i_1)_*\td_{m-1}(D_1 ) +(i_2)_*\td_{m-1} (D_2 )- \frac{T_{m-1}}{T_{m-2}}\cdot  (i_{12})_*\td_{m-2}(D_1\cap D_2 ),
\end{equation}
where by $i_1$, $i_2$, $i_{12}$, we denote the natural embeddings.

c) Suppose that $D\sim x-y$, with $x$ and $y$ smooth divisors on $W$.
Suppose also that there are smooth divisors $y_i$, $1\leq i\leq \delta$,
in the same linear equivalence class with $y$, such that, for each $k=1,\ldots , \delta$, the
scheme theoretic intersections $y_1\cap\cdots \cap y_k\cap y$,
$y_1\cap\cdots \cap y_k\cap x$ are smooth of pure codimension $k+1$. 
Denote by $i_x: x\hookrightarrow W$, $i_y: y\hookrightarrow W$,
$i_k: y_1\cap\cdots \cap y_k\cap x\hookrightarrow W$, $i'_k: y_1\cap\cdots \cap y_k\cap y\hookrightarrow W$,
the natural embeddings. Then,
we have
\begin{multline}
\ \ \td_{m}(D; W ) =(i_x)_*(\td_{m-1} (x ))-(i_y)_*(\td_{m-1}(y ))+\\
\ \ \ +\sum_{k=1}^{m-1}\frac{T_{m-1}}{T_{m-1-k}}\cdot \left[(i_k)_*  (\td_{m-1-k} (y_1\cap\cdots \cap y_k\cap x ))-
(i'_k)_*(\td_{m-1-k}  (y_1\cap\cdots \cap y_k\cap y ))\right] \ \ \ \
\end{multline} in $\CH^{m}(W)$. 
\end{prop}
 
\begin{Proof}
a) Since both $D$ and $W$ are smooth, we have
$
[i^* T_W ]=[T_D]+[\O_D(D)]
$
 in the Grothendieck group $\Kr_0(D)$. Therefore, by 
 the Whitney sum formula, we obtain
\begin{equation}\label{whit}
c_j([ i^*T_W ])=c_j(T_D)+c_1(\O_D(D))\cdot c_{j-1}(T_D)\ .
\end{equation}
 %Denote by $g: D\to C$ the composition $g=h\cdot i$.

For any polynomial $P$ (with integral coefficients) in the Chern classes,
a locally free coherent $\O_W$-sheaf $\F$ on $W$, and $n\geq 1$, we have
\begin{equation}
[D]^n\cdot P(\F)=i_*\left(c_1(\O_D(D))^{n-1}\cdot P(i^*\F)\right) .
\end{equation}
(On the right hand side, the Chern classes and the intersection are in $\CH^*(D)$.)
This identity implies
\begin{multline}\label{eqsmooth1}
\ \ \ Q_{ m}(c_1(T_W),   \ldots , c_{m-1}(T_W),    [D])=\\
=i_*\left(\left(T_{m-1}\cdot \left\{\frac{1-e^{-x}}{x}
\cdot \Td \right\}_{ m-1}\right)( c_1(i^*T_W), \ldots , c_{m-1}(i^*T_W),   c_1(\O_D(D))\right),
\end{multline}
where in the last expression the Chern classes are for bundles on $D$. 
 The usual expression of the Todd power series 
in terms of the Chern roots  gives
that the polynomial
$$
\left\{\frac{1-e^{-x}}{x} 
\cdot \Td\right\}_{ m-1}\ \ \in \ \Q[ c_1,  \ldots , c_{m-1},  x]
$$
is send to  
$
 \Td _{ m-1}( c_1, \ldots, c_{m-1})$ in $\Q[ c_1, \ldots , c_{m-1}]$
under the substitution $c_i\mapsto c_i+x\cdot c_{i-1}$,  $x\mapsto x$. This 
fact, together with (\ref{whit}), implies 
that the expression in (\ref{eqsmooth1}) above
 equals
\begin{equation*}
i_*\left(\{T_{m-1}\cdot \Td_{ m-1}\}( c_1(T_D),\ldots, c_{m-1}(T_D))\right).
\end{equation*}
This shows our claim.

b) Consider the identity
\begin{equation}\label{exp1}
1-e^{-a-b}=(1-e^{-a})+(1-e^{-b})-(1-e^{-a})(1-e^{-b}).
\end{equation}
Our claim  follows from  the identity (\ref{exp1}) and Lemma \ref{factorials}
by applying an argument similar to the proof of part (a).
(Under our assumptions, $D_1\cap D_2$ is a smooth divisor in $D_2$.)

c) Consider the formal identity
\begin{equation}\label{exp2}
1-e^b=-\sum_{j=1}^{\infty}(1-e^{-b})^j\ .
\end{equation}
This together with (\ref{exp1}) gives 
\begin{equation}\label{exp3}
1-e^{-a+b}=(1-e^{-a})-(1-e^{-b})+\sum_{j=1}^\infty((1-e^{-a})(1-e^{-b})^j- (1-e^{-b})^{j+1}).
\end{equation}
Our claim again follows using (\ref{exp2}) 
by  arguments as in the proofs of (a) and (b) above.\endproof
\end{Proof}

 \smallskip

\begin{prop}\label{vample2}
a) Under the assumptions  of   Proposition \ref{vample1} (b),
we have
\begin{equation*}
 [\O_W]- [\O_W(-D)]= [\O_{D_1}] + [\O_{D_2}]-   [\O_{D_1\cap D_2}].
\end{equation*}
b) Under the assumptions  of   Proposition \ref{vample1} (c), we have
\begin{equation*}
 [\O_W]- [\O_W(-D)]= [\O_x]- [\O_y]
+\sum_{k=1}^{\delta} \left( [\O_{y_1\cap\cdots \cap y_k\cap x}]-
 [\O_{y_1\cap\cdots \cap y_k\cap y}]\right)
\end{equation*}
in ${\rm K}_0(W)$. 

(For simplicity, we omit
from the notation the push forward along closed immersions
and write for example $\O_{D_1}$, $\O_x$, instead of $(i_1)_*(\O_{D_1})$,
$(i_x)_*(\O_x)$. )
\end{prop}

\begin{Proof}
For simplicity, we write $\O=\O_W$.

a) Under our assumptions, the restriction of $\O(-D_1)$ 
to $D_2$ is $\O_{D_2}(-D_1\cap D_2)$.
The result follows using the exact sequences
\begin{eqnarray*}
0\to \O(-D_1-D_2)\to \O(-D_1)\to \O_{D_2}(-D_1\cap D_2)\to 0  , \\
0\to \O_{D_2}(-D_1\cap D_2)\to \O_{D_2}\to \O_{D_1\cap D_2}\to 0  .\ \ \ 
\end{eqnarray*}

b)
The  exact sequence
\begin{equation*}
0\to \O(-x+y) \to \O(y)\to \O_x(y)\to 0
\end{equation*}
gives
\begin{equation}\label{deq1}
[\O]-[\O(-D)]=[\O_x]+([\O]-[\O(y)])-([\O_x]-[\O_x(y)])  .
\end{equation}
The exact sequences
\begin{equation*}
 0\to \O\to \O(y)\to \O_y(y)\to 0, \quad 0\to \O_y\to \O_y(y_1)=\O_y(y)\to \O_{y\cap y_1}(y)\to 0 ,
\end{equation*}
give
\begin{equation*}
[\O]-[\O(y)]=-[\O_y]-[\O_{y_1\cap y}(y)] .
\end{equation*}
Inductively, we now obtain
\begin{equation*}
 [\O]-[\O(y)]=- [\O_y]-\sum_{k=1}^{\delta} [\O_{y_1\cap\cdots \cap y_k\cap y}]  .
\end{equation*}
(Since  $y_1\cap\cdots \cap y_{\delta}\cap y$ is $0$-dimensional,  
$ [\O_{y_1\cap\cdots \cap y_{\delta}\cap y}(y)]=[\O_{y_1\cap\cdots \cap y_{\delta}\cap y}]$.)
The same argument also shows
\begin{equation*}
 [\O_x]- [\O_x(y)]= -\sum_{k=1}^{\delta} [\O_{y_1\cap\cdots \cap y_k\cap x}]  .
\end{equation*}
The  last two equations, combined with (\ref{deq1}), allow us to conclude the proof.\endproof
\end{Proof}
\medskip

Bertini's theorem implies that we can always 
satisfy the assumptions of Propositions \ref{vample1} (b) and \ref{vample2} (c):

\begin{prop}\label{Bertini}
For any Weil  divisor $D$ on $W$, we can  
find $x$, $y$, very ample smooth divisors on $W$ 
such that $D\sim x-y$. Given such $x$ and $y$
we can find in addition very ample smooth divisors  $y_1$, $\ldots $, $y_{\delta}$
with $y_i\sim y$ 
which are such that $y_1\cap\cdots \cap y_k\cap y$,
$y_1\cap\cdots \cap y_k\cap x$, are smooth of (pure) codimension $k+1$,
for all $k=1,\ldots , \delta$.\endproof
\end{prop}

\subsection{} Now suppose $D$ and $D'$ are two (arbitrary) Weil divisors 
on $W$. Apply Proposition \ref{Bertini} to $D$ and $D'$. We can write $D\sim X-Y$, $D'\sim X'-Y'$ with $X$, $Y$, $X'$, $Y'$ smooth very ample and find $Y_i\sim Y$, $Y_i'\sim Y'$, $k=1,\ldots ,\delta$, with the properties stated above. In addition, we
can arrange so that $X\cap X'$, $Y\cap Y'$ are both smooth and 
of pure codimension $2$. Since $X+X'$ and $Y+Y'$ are also very ample,
we can represent them by smooth divisors $U\sim X+X'$, $V\sim Y+Y'$.
Proposition \ref{vample1} applied to $D$, $D'$ and $D+D'$,
now gives
\begin{eqnarray*}
\td_m(D;W)&=&(i_X)_*(\td_{m-1}(X))-(i_Y)_*(\td_{m-1}(Y))+A,\\
\td_m(D';W)&=&(i_{X'})_*(\td_{m-1}(X'))-(i_{Y'})_*(\td_{m-1}(Y'))+A',\\
\td_m(D+D';W)&=&(i_U)_*(\td_{m-1}(U))-(i_V)_*(\td_{m-1}(V))+A'' .
\end{eqnarray*}
Similarly, Proposition \ref{vample1} (a) and (b) applied to $U\sim X+X'$, $V\sim Y+Y'$,
gives
\begin{eqnarray*}
(i_U)_*(\td_{m-1}(U))=\td_m(U;W)&=&(i_X)_*(\td_{m-1}(X))+(i_{X'})_*(\td_{m-1}(X'))+B ,\\
(i_Y)_*(\td_{m-1}(V))=\td_m(V;W)&=&(i_Y)_*(\td_{m-1}(Y))+(i_{Y'})_*(\td_{m-1}(Y'))+B' .
\end{eqnarray*}
Here $A$, $A'$, $A''$ and $B$, $B'$ 
are integral linear combinations of   
classes of the form
$$
\frac{T_{m-1}}{T_{m-1-l}}\cdot (i_Z)_*(\td_{m-1-l}(Z))
$$
with $i_Z: Z\hookrightarrow W$   a smooth subvariety of 
codimension $l+1\geq 2$ in $W$. 
By combining  the corresponding equations,
we obtain 
\begin{multline}\label{320}
\td_{m}(D+D';W ) =\\ =\td_{m}(D;W ) +\td_{m}(D';W )+\sum_{l=1}^{m-1} \left[\sum_{Z} a_{Z}\frac{T_{m-1}}{T_{m-1-l}}\cdot (i_Z)_*\td_{m-l-1}(Z )\right]
\end{multline}
in $\CH^m(W)$. The sum in the bracket is over the set of  smooth   subvarieties $i_Z: Z\hookrightarrow W$ of codimension $l+1$
and $a_Z$ are integers which are almost always $0$. 

Similarly, the same argument applied to the  equations 
obtained by Proposition \ref{vample2} gives
\begin{equation}\label{321}
[\O]-[\O(-D-D')]%=\ \ \ \ \ \ \ \ \ \ \ \ \ \ \ \ \ \ \ \ \\
=([\O]-[\O(-D)])+([\O]-[\O(-D')])+\sum_{l=1}^{\delta} \left( \sum_{Z}b_{Z}\cdot  [\O_Z]\right)
\end{equation}
in ${\rm G}_0(W)={\rm K}_0(W)$.
In fact, the parallel expressions in the statements
of Propositions \ref{vample1}, \ref{vample2}, 
allow us to observe that if $2\leq {\rm codim}(Z)\leq m$, then $ b_Z=a_Z$. 
%If $D_1$, $D_2$ are effective
%and irreducible we have
%\begin{equation}
% [\O]- [\O(-D_1-D_2)]= [\O_{D_1}]+ [\O_{D_2}]+\sum_{l=1}^{\delta} \left(\sum_{Z}  b_{Z}\cdot  [\O_Z]\right)\ .
%\end{equation}

\bigskip

\section{Projective bundles, blow-ups and embeddings}\label{bundles}

\setcounter{equation}{0}

In this section, we show that integral Riemann-Roch holds 
for  projective bundles, for closed immersions and for blow ups
along smooth centers. The proofs mostly follow  the 
standard  arguments of ``Riemann-Roch algebra" ([FL]);
essentially, we will observe that the integrality of the expressions
involved is preserved. 

\subsection{} \label{projbundles}
Suppose that $\E$ is a locally free coherent sheaf of rank $r+1$ over the smooth quasi-projective variety
$Y$ and denote by $p: {\bf P}(\E)={\bf Proj}({\rm Sym}(\E))\to Y$ the corresponding projective bundle.
%The dimension of ${\bf P}(\E)$ is $r+m$.
Denote by $ \O_{{\bf P}(\E)}(1)$  the Serre invertible sheaf
and set, as usual,  $ \O_{{\bf P}(\E)}(a)= \O_{{\bf P}(\E)}(1)^{\otimes a}$.
Recall ([FL, V, Theorem 2.3]) that the Grothendieck ring 
${\rm K}_0({\bf P}(\E))$ is isomorphic to
\begin{equation*}
{\rm K}_0(Y)[T]/(T^{r+1}-[\E]\cdot T^{r}+ \cdots +(-1)^{r+1}[\wedge^{r+1}(\E)])\, ,
\end{equation*}
with $l=[\O_{{\bf P}(\E)}(1)]$  
corresponding to the class of the element $T$ in the quotient.
Under this isomorphism, the pull-back $p^*:{\rm K}_0(Y)\to {\rm K}_0({\bf P}(\E))$ is identified
with $a\mapsto a\cdot T^0$.
Similarly,   for the Chow ring we have
([F1, Theorem 3.3]):
\begin{equation}\label{chowproj}
\CH^*({\bf P}(\E))\simeq \CH^*(Y)[T]/(T^{r+1}-c_1(\E)T^{r}+ \cdots +(-1)^{r+1}c_{r+1}(\E))\, .
\end{equation}
Under this isomorphism, the grading of $\CH^*({\bf P}(\E))$ 
corresponds to the grading given by setting $\deg(a\cdot T^i)=\deg(a)+i$,
for $0\leq i\leq r$.
(Here, the class of $T$  
corresponds to the first Chern class of $\O_{{\bf P}(\E)}(1)$.)
The pull-back $p^*:\CH^*(Y)\to \CH^*({\bf P}(\E))$ is identified
with $a\mapsto a\cdot T^0$.
By [F1, Prop. 3.1 (a)], the push 
forward map $p_*: \CH^*({\bf P}(\E))\to \CH^{*-r}(Y)$ is identified,  
 under (\ref{chowproj}), with
$p_*(a\cdot T^r)=a$, $p_*(a\cdot T^i)=0$ if $0\leq i\leq r-1$.

\begin{prop}\label{bun1}
We have $p_*[\O_{{\bf P}(\E)}]=[\O_Y]$, and $p_*[\O_{{\bf P}(\E)}(a)]=0$
if $-r\leq a <0$.
\end{prop}

\begin{Proof}
This is well-known; see for example [FL, V \S 2].\endproof
\end{Proof}
\smallskip

\begin{thm}\label{projcor}
Theorem \ref{main} is true 
for the projective bundle $p: {\bf P}(\E)\to Y$,
i.e
\begin{equation}\label{projeq}
\frac{T_{n+r}}{T_n}\cdot \ctd_{n}(p_*[\F], Y)=
p_*(\ctd_{n+r}( \F , {\bf P}(\E)))
\end{equation}
for all $n\geq 0$. 
\end{thm}

\begin{Proof}
This closely follows the proof of Riemann-Roch for elementary projections given in [FL,
II \S 2]:
The description  above implies that  ${\rm K}_0({\bf P}(\E))$
is generated as a ${\rm K}_0(Y)$-module by the classes
$l^a=[\O_{{\bf P}(\E)}(a)]$ for $a=-r,-r+1,\ldots, -1, 0$.
Using Proposition \ref{multiply}, we see that it is enough to 
show the identity (\ref{projeq}) for these classes.
 By Proposition \ref{bun1}, it is enough to show 
\begin{equation}
p_*(\ctd_{n+r}(l^a, {\bf P}(\E)))=
\begin{cases}
0,&  \hbox{\rm if  } -r\leq a<0\\
\displaystyle{\frac{T_{n+r}}{T_n}}\cdot \td_n(T_Y),&  \hbox{\rm if   } a=0.
\end{cases}
\end{equation}
 
 There is an exact sequence  
 ([FL, IV Prop. 3.13])
\begin{equation}\label{t2}
0\to \O_{{\bf P}(\E)}\to p^*\E^{\vee}\otimes \O_{{\bf P}(\E)}(1) 
\to T_{{\bf P}(\E)}\to p^*T_{Y} \to 0\ .
\end{equation}
 \begin{comment}
There is an exact sequence of tangent bundles
\begin{equation}\label{t1}
0\to T_{{\bf P}(\E)/Y}\to T_{{\bf P}(\E)}\to p^*T_{Y}\to 0\ 
\end{equation}
and an exact sequence ([FL, IV Prop. 3.13])
\begin{equation}\label{t2}
0\to \O_{{\bf P}(\E)}\to p^*\E^{\vee}\otimes \O_{{\bf P}(\E)}(1)\to T_{{\bf P}(\E)/Y} \to 0\ 
\end{equation}
\end{comment}
For simplicity, we will denote $p^*\E^{\vee}\otimes \O_{{\bf P}(\E)}(1)$ by $\E^\vee(1)$.
The exact sequence  (\ref{t2}) implies that $c_i(T_{{\bf P}(\E)})=c_i(\E^\vee(1)+p^*T_Y)$ for all $i$. This, together with (\ref{mulTodd}) and the projection formula, proves that
it is enough to show
\begin{equation}
p_*\left[\sum_{j=0}^{n'}\frac{T_{n'}}{j!\cdot T_{n'-j}}\cdot \ch_j(l^a)\cdot \td_{n'-j}(\E^\vee(1))\right]=
\begin{cases}
\ 0,&  \hbox{\rm if  } -r\leq a<0,\, \hbox{\rm or }n'\neq r\\
{T_{r}}\cdot 1,&  \hbox{\rm if   } a=0,\, \hbox{\rm and } n'=r\,.
\end{cases}
\end{equation}
The  proof can now be completed as in [FL,
II \S 2]: Indeed, the term in the bracket above is  
the (integral) $n'$-th degree  term   
\begin{equation*}\label{Howe}
T_{n'}\cdot \left\{e^{aT}\cdot \prod_{i=1}^{r+1}\frac{T-x_i}{1-e^{-(T-x_i)}}\right\}_{n'}
\end{equation*}
of the (symmetric) formal power series 
$$
H_a(T, \{x_i\})=e^{aT}\cdot  \prod_{i=1}^{r+1}\frac{T-x_i}{1-e^{-(T-x_i)}}\in \Q[[x_1,\ldots ,x_{r+1}, T]]\, ,
$$ 
evaluated  at $T=c_1(l)$ and $x_i=$ Chern roots of $\E$.
By loc. cit. Lemma 2.3 and its proof, the relation 
$\prod_{i=0}^{r+1}(T-x_i)=0$ implies that we can write 
$$
H_a(T, \{x_i\})=\sum_{j=0}^r f_{j,a}(c_1,\ldots , c_{r+1})\cdot T^j ,
$$
with $f_{j,a}\in \Q[c_1,\ldots , c_{r+1}]$ and $f_{r,a}=1$
if $a=0$, $f_{r,a}=0$ if $-r\leq a<0$. (Here $c_i$ is the $i$-th elementary symmetric function
of  $x_1,\ldots ,x_{r+1}$.) The desired conclusion now follows using the above description of   $p_*: \CH^*({\bf P}(\E))\to \CH^{*-r}(Y)$.%  (cf.   proof of [FL, II, Theorem 2.2]).
 \endproof
\end{Proof}

\subsection{}\label{blow}
Let $b: \ti X\to X$ be the blow up of the smooth quasi-projective variety $X$
along the smooth subvariety $Z$. Denote by $i:Z\hookrightarrow X$ 
the embedding and by $N=N_{Z|X}$ the normal sheaf of $Z$ in $X$.

We now recall some aspects
of the construction of   ``deformation to the normal cone"
([F1, 5.1], or [FL, IV \S5]). %Our conclusions below are similar to [Fu2, Lemma 1].  
%For simplicity, we will set $N^\vee\oplus 1=N^\vee\oplus \O_Z$.
Consider the blow-up of $\pi: W\to X\times {\bf P}^1$ along $Z\times \{\infty\}$. Let   $D$ be the divisor on $W$
which is the preimage $D =\pi^{-1}(X\times\{0\})$ of $X\times \{0\}$; we can identify $D$ with $X$. 
The divisor $D$ is linearly equivalent on $W$  to the divisor $D'$ given by the preimage $\pi^{-1}(X\times\{\infty\})$; 
 $D'$ is the sum ${\ti X}+{\bf P}(N^\vee\oplus \O_Z) $ of two smooth irreducible components: 
the projective bundle ${\bf P}(N^\vee\oplus \O_Z)$
over $Z$ and the blow-up $\ti X$. The scheme theoretic intersection $\ti X\cap {\bf P}(N^\vee\oplus \O_Z)$
is the (smooth) projective bundle ${\bf P}(N^\vee)$ over $Z$ (this is the exceptional locus of the blow-up $\ti X\to X$).
%Assume now that $X$ supports a projective morphism $f: X\to S$ of relative dimension $d$.
%Denote by $g: W\to S$ the composition of $\pi$ with $pr_1\cdot f: X\times {\bf P}^1\to S$
%and by $\ti f: \ti X\to S$, $f_Z: Z\to S$ the morphisms obtained 
%by composing the blow-up morphism $b: \ti X\to X$ and $Z\hookrightarrow X$ with $f$.

Using Proposition \ref{vample1} (a) and (b), we obtain
\begin{multline}\label{433}
(i_X)_*(\td_{m}(X))=\td_{m+1}(D; W) = \td_{m+1}(D'; W)=\\
=(i_{\ti X})_* \td_{m} (\ti X) +(i_{{\bf P}(N^\vee\oplus \O_Z)} )_* \td_{m } ({\bf P}(N^\vee\oplus \O_Z) )-
\frac{T_{m}}{T_{m-1}}\cdot (i_{{\bf P}(N^\vee)})_* \td_{m-1} ({\bf P}(N^\vee))  
\end{multline}
in $\CH^{m+1}(W)$.
Consider the composition $q:=pr_1\cdot \pi: W\to X\times {\bf P}^1\to X$. 
Observe that $q\cdot i_X={id}_X$, $q\cdot i_{\ti X}=b$.
Also $q\cdot i_{{\bf P}(N^\vee\oplus \O_Z)}$, $q\cdot i_{{\bf P}(N^\vee )}$,
are the compositions of the projective bundles ${\bf P}(N^\vee\oplus \O_Z)\to Z$,
resp.
${\bf P}(N^\vee)\to Z$, with $i_Z:Z\hookrightarrow X$. 
Now apply the push-forward  homomorphism $q_*$ to the identity (\ref{433}).
Using Proposition \ref{bun1}, Theorem \ref{projcor} for $\F=$ the structure
sheaf,  and the above observations, we obtain
\begin{equation}\label{450}
\td_{m}(X )=b_*(\td_{m}(\ti X ))\ 
\end{equation}
in $\CH^{m}(X)$, for all $0\leq m\leq \dim(X)$. 
 
A similar argument, using Proposition \ref{vample2} (a) and Proposition \ref{bun1},  gives
\begin{equation}\label{445}
 [\O_{ X}] =  b_* [\O_{\ti X}]\, .
\end{equation}
(This also follows from [SGA6, VII, Prop. 3.6].)
 
\begin{comment}It follows from (\ref{444}) and {\ref{445}) 
that the main theorem (for $\F$ equal to the structure sheaf and $\ch_n$)
is true for $X\to S$  if and only if it is true for a blow-up $\ti X\to S$
along a smooth subvariety $Z\subset X$. 
\end{comment}

\subsection{} \label{immers}

Let $i: Z\hookrightarrow X$ be as above. Set $r=-d_i=\dim(X)-\dim(Z)$.

\begin{thm}\label{JoRR} Given $n\geq 0$, and $\F$  
a locally free coherent  $\O_Z$-sheaf on $Z$, we have
\begin{equation}\label{imeq}
\frac{T_{n}}{T_{n-r}}\cdot i_*(\ctd_{n-r}(\F, Z))=\ctd_{n}(i_*[\F], X)\,.
\end{equation}
\end{thm}
%(Of course, when $n<r$, the statement is that $\ctd_{n }(i_*[\F]; X)=0$.)
\begin{Proof}
This can be deduced from the Riemann-Roch theorem 
``without denominators" for regular immersions 
([J]). Here we give a  direct argument using
the technique of deformation to the normal cone. 
We have the  formal identity  
\begin{equation}\label{formal}
\sum_{j=0}^r (-1)^j\sum_{i_1<\cdots <i_j}e^{-x_{i_1}-\cdots -x_{i_j}}=\prod_{i=1}^r(1-e^{-x_i})
=x_1\cdots x_r\cdot \prod_{i=1}^r\frac{1-e^{-x_i}}{x_i}\,.
\end{equation}
The degree $m$ part of the left hand side
is zero when $m<r$ and equal to
$$
x_1\cdots x_r\cdot \left\{\prod_{i=1}^r\frac{1-e^{-x_i}}{x_i} \right\}_{m-r}
$$
if $m\geq r$. The denominator of the degree $m$ part 
of the left hand side of (\ref{formal})
divides $m!$ and so
$$
\td^{inv}_{m-r}:=m!\cdot \left\{\prod_{i=1}^r\frac{1-e^{-x_i}}{x_i} \right\}_{m-r}
$$
is a symmetric homogeneous polynomial with integral coefficients.
As such, it can be expressed as an integral polynomial in the elementary symmetric functions $c_1, \ldots , c_r$ of the variables $x_1,\ldots ,x_r$.
Denote by $\td^{inv}_{m-r}(\Gg)$ the result of evaluating 
$\td^{inv}_{m-r}$ at the Chern classes $c_i=c_i(\Gg)$. 

\begin{prop}\label{prop44} Recall that $N$ is the normal sheaf of $Z$ in $X$.  We have
\begin{equation}\label{chclaim}
\ch_m(i_*[\F])= {\sum_{l=r}^m \frac{m!}{(m-l)!\cdot l!}\cdot  i_*(\ch_{m-l}(\F)\cdot \td^{inv}_{l-r}(N))} 
\end{equation}
in $\CH^m(X)$, for all $m\geq r$, while $\ch_m(i_*[\F])=0$ for $m<r$.
\end{prop}

\begin{Proof} This follows the proof of Riemann-Roch 
theorem for regular immersions in [F1, \S 15.2] mutatis-mutandis.
As in loc. cit. p. 287-288, we see  that
it is enough to prove the statement in the model situation in which
$X={\bf P}(N^\vee\oplus \O_Z)$ and $i: Z\hookrightarrow {\bf V}(N)={\bf Spec}({\rm Sym}(N^\vee))\subset {\bf P}(N^\vee\oplus \O_Z)=X$,
where $Z\hookrightarrow {\bf V}(N)$ is the zero section of the bundle ${\bf V}(N)\to Z$. The argument for showing (\ref{chclaim}) in this case, is  similar
to the argument in loc. cit., p. 282-283:
%Since in our case we have to take care not to ``destroy integrality",
%we repeat some of the arguments here. 
Let $p: X={\bf P}(N^\vee\oplus \O_Z)\to Z$ be the projection, 
let $Q$ be the universal kernel sheaf on ${\bf P}(N^\vee\oplus \O_Z)$
and let $s$ be the section of ${\bf V}(Q^\vee)$ determined by the projection of $Q$
to the trivial factor of $p^*(N^\vee\oplus \O_Z)$. The (scheme) theoretic zero locus of $s$ is $Z$. Since we have $\F\simeq i^*(p^*\F)$,  
by using the projection formula and (\ref{mulChern}), we see that it is enough to show the claim for 
$\F=\O_Z$. The coherent sheaf
$i_*\O_Z$ on $X$ has a Koszul resolution
\begin{equation}
0\to \wedge^r Q  \to \cdots \to Q    \xrightarrow{s^*} \O_X\to i_*\O_Z\to 0\ .
\end{equation}
This, together with the identity (\ref{formal}), now gives
\begin{equation}\label{im1}
\ch_{m }(i_*[\O_Z])
=   \sum_{j=0}^{r} (-1)^j\ch_m(\wedge^j Q) ={\begin{cases}
  c_{r}(Q)\cdot  \td^{inv}_{m-r}(Q), \ \hbox{\rm if\ }  m\geq r, \\
  0, \ \ \ \ \ \ \ \ \ \ \ \ \ \ \ \ \ \ \ \ \ \  \ \hbox{\rm if\ }  m<r .
\end{cases}}
\end{equation}
Observe that $c_r(Q)$ is represented by the zero locus  $Z$ 
of the section $s: X\to {\bf V}(Q^\vee)$. Therefore, we have
\begin{equation}\label{im4}
c_{r}(Q)\cdot \td_{m-r}^{inv}(Q)=i_*(\td^{inv}_{m-r}(i^*Q))
\end{equation}
in $\CH^m(X)$.
%, where on the right hand side
%$\td^{inv}_{m-r}(i^* Q )\in \CH^{m-r}(Z)$ is formed using 
%Chern classes on $Z$. 
By our construction,  $i^*Q$ is isomorphic to $N$ and our claim   for $\F=\O_Z$ and $i: Z\to X={\bf P}(N^\vee\oplus\O_Z)$  follows from (\ref{im1}) and (\ref{im4}). By the above discussion, the 
proof of the proposition also follows.\endproof
\end{Proof}
\smallskip

Let us now complete the proof of Theorem \ref{JoRR}. Observe that 
$[N]+[T_Z]=[i^*T_X]$ in ${\rm K}_0(Z)$. This, together with the definition of $\td_{l}^{inv}(N)$ via the inverse of the Todd power series, implies
\begin{equation}\label{im3}
\frac{T_{m}}{T_{m-r}}\cdot \td_{m-r}(T_Z)=  \sum_{j=0}^{m-r} \frac{T_{m }}{(j+r)!\cdot  T_{m-r-j}}\cdot \td_{j}^{inv}(N)\cdot \td_{m-r-j}(i^*T_X)\, . 
\end{equation}
The identity (\ref{imeq}) of Theorem \ref{JoRR} now follows from Proposition \ref{prop44} 
by using (\ref{im3}) and the projection formula.
\endproof \end{Proof}

\subsection{} \label{moreblow}

We continue with the assumptions and   notations of the previous paragraphs. 
In particular, $b: \ti X\to X$ is the blow-up of $X$
along $Z$.  

\begin{thm}\label{blowRR}
Let $\F$ be a locally free coherent $\O_{\ti X}$-sheaf on   $\ti X$.
Then for $n\geq 0$ we have
\begin{equation}\label{morebloweq}
b_*(\ctd_{n}(\F, \ti X))=\ctd_n(b_*[\F], X)
\end{equation}
in $\CH^n(X)$.
\end{thm}

\begin{Proof}
We have the (commutative) blow-up diagram
\begin{equation}
\begin{matrix}
\ti Z={\bf P}(N^\vee)\ \ \ &\xrightarrow{\ j\ } & \ti X\\
q \downarrow\ \  && \ \ \ \downarrow b\ \\
Z&\xrightarrow{\ i\ } & \ X\,.
\end{matrix}
\end{equation}
It   follows from (\ref{450}), (\ref{445}) and Proposition \ref{multiply} that (\ref{morebloweq}) is true for $\F$ of the form $b^*\Gg$, with 
$\Gg$ a locally free coherent $\O_X$-sheaf on $X$.
By [SGA6, VII, Th. 3.7], each element of ${\rm K}_0(\ti X)={\rm G}_0(\ti X)$ can be written in the form $a=b^*(a')+j_*(z)$ with $a'\in {\rm K}_0(X)$, 
$z\in {\rm K}_0(\ti Z)$. Since both sides of (\ref{morebloweq})
are additive, it remains to prove the equality (\ref{morebloweq}) for $\F=j_*\Hh$, where $\Hh$ is a 
locally free coherent $\O_{\ti Z}$-sheaf on $\ti Z$.
Using $b_*\cdot j_*=(b\cdot j)_*=(i\cdot q)_*=  i_*\cdot q_*$,   Theorem \ref{JoRR} for the embeddings $j$ and $i$,
and Theorem \ref{projcor} for $q: \ti Z={\bf P}(N^\vee)\to Z$,  we
obtain
\begin{multline*}
b_*(\ctd_{n}(j_* \Hh , \ti X))
 = b_*\left( \frac{T_{n}}{T_{n-1}}\cdot
j_*(\ctd_{n-1}(\Hh, \ti Z))\right)=  \\
  =\frac{T_n}{T_{n-1}}\cdot i_*\left(  q_*(\ctd_{n-1}(\Hh, \ti Z))\right) =\ \ 
   \ \ \ \ \ \ \ \ \ \ \ \ \ \ \ \ \ \ \ \\
 \ \ \   \ = \frac{T_n}{T_{n-1}}\cdot   i_*\left(\frac{T_{n-1}}{T_{n-1-(r-1)}}\cdot\ctd_{n-1-(r-1)}(q_*[\Hh],  Z)\right)=
 \\ =\frac{T_n}{T_{n-r}}\cdot   i_*\left( \ctd_{n-r}(q_*[\Hh],  Z)\right)=\ \ \ \  \ \ \ \ \ \ \ \ \ \ \ \ \ \ \ \ \ \\  =\ctd_{n }(i_*q_*[\Hh],  X)= 
 \ctd_{n }(b_*[j_* \Hh ],  X).   \ \ \ \ \ \ \ \ \ \  \ \ \ \ \ \ \ \ \ \ \ \ \ \ \ \ \ \ \ \ \ \ \ \ \ \ \ 
\end{multline*}\endproof
\end{Proof}

\section{Factorization}\label{factor}

\setcounter{equation}{0}

Let $f: Y\to S$ and $f': Y'\to S$ be projective morphisms between 
smooth quasi-projective varieties.
We will say that  
{\sl $f$, $f'$ are birationally isomorphic over 
$S$} if the following is true:
We have $f(Y)=f'(Y')$ and there is an isomorphism of the function fields 
$a: k(Y')\xrightarrow{\sim} k(Y)$
which commutes with the $k(f(Y))=k(f'(Y'))$-algebra structures given by $f'$, $f$.

Now assume that $f$, $f'$ are birationally isomorphic over 
$S$ and let us write $\phi: Y-\to Y'$ for the corresponding birational map. 
Let $U\subset Y$ be the largest open subscheme of $Y$ such that $\phi_{|U}: U\to \phi(U)$
is an isomorphism; in what follows, we will implicitly identify $U$ and $\phi(U)$.

\begin{thm}\label{factorThm}
([AKMW]) There is a finite sequence of
birational maps between smooth quasi-projective varieties
\begin{equation*}
Y=Y_0-\xrightarrow{\phi_1} Y_1-\xrightarrow{\phi_2}  \cdots -\xrightarrow{\phi_{n-1}}Y_{n-1}
-\xrightarrow{\phi_{n}} Y_n=Y'
\end{equation*}
over $k$ such that:

\begin{enumerate}
\item $\phi=\phi_n\cdot\cdots \cdot \phi_2\cdot \phi_1$,

\item the $\phi_i$'s are isomorphisms on $U$,

\item for each $i$, either $\phi_i: Y_{i-1}-\to Y_i$ or $\phi^{-1}_i: Y_i-\to Y_{i-1}$ is obtained by blowing up
a nonsingular subscheme disjoint from $U$,

\item there is an index $i_0$ such that for all $i\leq i_0$, resp. $i\geq i_0$, the birational map
$\phi(i):=\phi^{-1}_1\cdots \cdot \phi^{-1}_{2} \cdot \phi^{-1}_i: Y_i-\to Y_0=Y$, resp.
$\phi(i):=\phi_i\cdot \phi_{i+1}\cdot \cdots \cdot \phi_n: Y_i-\to Y_n=Y'$, is a projective morphism,

\item the varieties $Y_i$ support  projective morphisms $f_i: Y_i\to S$
such that, for each $i$, the blow-up morphism $\phi_i: Y_{i-1}-\to Y_i$ or $\phi^{-1}_i: Y_i-\to Y_{i-1}$ given by item (c) is  
a 
morphism over $S$  (i.e commutes with $f_i$ and $f_{i-1}$). 
\end{enumerate}
\end{thm}

\begin{Proof}
This follows from the ``functorial weak factorization theorem"
(Theorem 0.3.1, cf. Remark (1)) of [AKMW]. The result in loc. cit.
gives varieties $Y_i$ and birational maps $\phi_i$ with the properties (a)-(d).
It remains to observe that item (e) follows from (a)-(d). Indeed, we can
 define an $S$-structure on $Y_i$ as follows, using (d): If 
 $i\leq i_0$, then we set $f_i=f\cdot \phi(i)$. If $i\geq i_0$,
 then we set $f_i=f'\cdot \phi(i)$. We can see that these $f_i$ satisfy the requirements
 of item (e).\endproof
\end{Proof}

\begin{cor}\label{corfactor}
If $f$, $f'$ are birationally isomorphic over 
$S$, then integral Riemann-Roch holds for $(f, n)$ if and 
only if it holds for $(f', n)$.
\end{cor}

\begin{Proof}
Using Theorem \ref{factorThm}, we see that it is enough 
to show the statement under the additional assumption 
that the birational map $\phi: Y'\to Y$ is obtained by blowing-up
$Y$ along a smooth center: By Theorem   \ref{blowRR},
integral Riemann-Roch holds for   $\phi$.
If integral Riemann-Roch holds for $(f, n)$ then,
by Proposition \ref{compose}, integral Riemann-Roch holds for $( f\cdot \phi, n)=(f', n)$. To show the converse, let
$\F$ be a locally free coherent $\O_{Y}$-sheaf on $Y$.
We will show that the integral Riemann-Roch identity for $\phi^*\F$ and $(f', n)$
implies the integral Riemann-Roch identity  
for $\F$ and $(f, n)$.
Suppose first   $d=d_f=d_{f'}\geq 0$. Since $f'_*=f_*\cdot \phi_*$, we obtain:  
\begin{equation*}
\frac{T_{n+d}}{T_n}\cdot 
\ctd_{n}(f_*\phi_* [\phi^*\F], S)=f_*\phi_*(\ctd_{n+d}(\phi^*\F, Y')) .
\end{equation*}
However, we have $\phi_*[\phi^*\F]=[\F]$ by [SGA6, VII Prop. 3.6.],  
while $\phi_*(\ctd_{n+d}(\phi^*\F, Y'))= \ctd_{n+d}(\phi_*[\phi^*\F],  Y)$
by Theorem \ref{blowRR}. Therefore, we obtain that integral Riemann-Roch 
holds  for $\F$ and $(f, n)$. The argument for $d<0$ is similar.\endproof 
\end{Proof}

\section{Completion of the proof}\label{nondon}
\setcounter{equation}{0}

We will now complete the proof of  Theorem \ref{main}
by induction on the degree $n$.

When $n=0$, since $\CH^0(S)=\Z$ is torsion-free, Theorem \ref{main} follows  from the standard Grothendieck-Riemann-Roch 
theorem. (More directly, we can deduce this case from Theorems \ref{projcor} and \ref{JoRR}
as in Grothendieck's proof.)
We  suppose that Theorem \ref{main}  %(i.e the equation (\ref{RReqA}) for $d_f\geq 0$and (\ref{RReqB}) for $d_f<0$) 
is known for all $f$, and for 
all degrees $n'<n$.

\subsection{}\label{noDom} Let $f: X\to S$ be a {\sl non-dominant} projective morphism between
the smooth quasi-projective varieties $X$ and $S$. We will
show that integral Riemann-Roch for $(f, n)$ follows 
from our inductive hypothesis.

Denote by $T$ the image $f(X)$; this is a subvariety of $S$ of codimension $r>0$.
By embedded resolution of singularities ([H]), 
we can find a projective birational  morphism $b: \ti S\to S$, which is obtained 
as a succession of blow-ups along smooth centers, so that the  strict transform 
$\ti T\subset \ti S$ of $T$  is smooth. Consider now the base change $X\times_T\ti T$
of $X\to T$ along $b: \ti T\to T$. By applying resolution of singularities to $X\times_T\ti T$, we can find a projective birational morphism $\phi: \ti X\to X\times_T\ti T\to X$ such that $\ti X$ is  smooth. Thus, we obtain a commutative diagram
\begin{equation*}
\begin{matrix}
 \ti X&\xrightarrow{\ti g}  & \ti T &\xrightarrow {\ti i} &\ti S\\
\phi\downarrow\ \ & &\ \ \downarrow b && \ \ \downarrow b \\
  X&\xrightarrow{  g}  &  T &\xrightarrow { i}  & S\\    
\end{matrix}
\end{equation*}
with $i\cdot g=f$. Set $\ti f=\ti i\cdot \ti g$. 
By Theorem \ref{JoRR}, integral Riemann-Roch holds for $(\ti i, n)$. 
By our induction hypothesis, it also holds for $(\ti g, n-r)$. Proposition \ref{compose} applied to the composition
$\ti f=\ti i\cdot \ti g$ allows us to conclude 
 that integral Riemann-Roch holds for $(\ti f, n)$.
Theorem \ref{blowRR} and Proposition \ref{compose}
applied to the composition $b\cdot \ti f$ now imply that integral Riemann-Roch
also holds for $(b\cdot \ti f, n)$. 
Since   $f\cdot \phi=b\cdot \ti f$, the morphisms $b\cdot \ti f: \ti X\to S$ and $f: X\to S$ are birationally isomorphic over $S$ via $\phi$.
Corollary \ref{corfactor} now shows that integral Riemann-Roch for $(f, n)$ follows.  

\subsection{}\label{6b}

We will now deduce  integral Riemann-Roch for $(f, n)$ 
from our inductive hypothesis on $n$, by using induction on the relative dimension
$d=d_f$. If $d<0$, then 
$f$ is not dominant and the result follows by the previous paragraph.
Hence, we  can assume that $d\geq 0$. 
%If $n=d=0$ then the statement was also shown above.  
We now suppose that 
integral Riemann-Roch for $(g, n)$ holds for all 
$g$ of relative dimension $d_g<d$; recall that we also assume that integral Riemann-Roch
for $n'$ holds  for  all morphisms, provided that $n'<n$. 

For a fixed smooth quasi-projective variety $S$
over $k$, let $\SS_{d}$ be the set of $S$-isomorphism classes $[f]$ of projective $k$-morphisms $f: X\to S$ 
of relative dimension $d\geq 0$ with $X$ smooth over $k$.
We will first show that integral Riemann-Roch holds for $(f, n)$ and the structure 
sheaf $\F=\O_X$, for all such $f$. In other words, we will show that
the ``error" function 
$E : \SS_d\to \CH^n(S)$ given by
\begin{equation}\label{errf}
 E ([f:X\to S])=  \frac{T_{d+n}}{T_n}\cdot \ctd_{n}(f_*[\O_X], S)-f_*(\td_{d+n}(X))\, 
\end{equation}
vanishes.
The induction hypothesis and our results in the previous sections give:
\medskip
 
(i) If $f: X\to S$ is $S$-isomorphic to a composition of the form $X={\bf P}(\E)\to Y\to S$,
with ${\bf P}(\E)\to Y$ a projective bundle   and $\dim(Y)<\dim(X)$,
then  $ E ([f:X\to S])=0$. (This follows from the induction hypothesis,
Proposition \ref{bun1} and
Theorem \ref{projcor}.) 
\smallskip

(ii) If $Z$ is a smooth subvariety of $X$, 
$$
E ([f: X\to S])=E ([\ti f: \ti X\to S]),
$$
with $\ti X$ the blow-up of $X$ along $Z$.
(This follows from (\ref{450}) and (\ref{445}).) 
\smallskip

In fact, Theorem \ref{factorThm} allows us to strengthen (ii)
(cf. proof of Corollary \ref{corfactor}):
 \smallskip
 
(ii)$'$ If $f: X\to S$ and  $f': X'\to S$ are birationally isomorphic over $S$ then
$$
 E ([f: X\to S])=E ([f': X'\to S]) .
 $$ 

(iii) Suppose that $g: W\to S$ is in $\SS_{d+1}$ and restrict $E$ to the smooth 
codimension $1$ closed subschemes of $W$, i.e to $f=g\cdot i$ with
$i: X\hookrightarrow W$ of codimension $1$. Proposition \ref{vample1} (a)
implies that
we can extend $E$ to Weil divisors by defining 
\begin{equation}\label{Eg}
E^g(D)=\frac{T_{d+n}}{T_n}\cdot\ctd_{n}(g_*([\O_W]-[\O_W(-D)]), S)-g_*(\td_{d+1+n}(D; W )).
\end{equation}
We can see that $E^g$ respects linear equivalence.
Using (\ref{320}) and (\ref{321}) for $m=d+n+1$ (including the fact that $b_Z=a_Z$) and our induction hypothesis, we see that
$$
E^g:  {\rm Pic}(W)\to \CH^n(S) 
$$
is a group homomorphism. 
\smallskip

(iv) If $f$ is {\sl not} dominant
then  
$E([f: X\to S])=0$. 
 
\subsection{}\label{last} We continue with our inductive proof.
The argument here was inspired by the calculation of the cobordism ring of a point by Levine and Morel [LM, Theorem 4.3.7].
Let $f: X\to S$ be a   projective morphism
as above which is dominant.
Then the generic fiber $X_K$ of $f$ is a smooth projective variety over $K=k(S)$
of dimension $d$.  Using the primitive element theorem (recall that ${\rm char}(k)=0$), we see that $X_K$ is birationally
isomorphic to a closed irreducible hypersurface $Y\subset {\bf P}^{d+1}_K$.
Denote by $\hat Y$ the Zariski closure of $Y$ in ${\bf P}^{d+1}_S={\bf P}^{d+1}_k\times_k S$; this is a divisor in ${\bf P}^{d+1}_S$ and affords a 
dominant projective morphism $\hat Y\to S$. For simplicity of notation,
we set $P={\bf P}^{d+1}_k\times_k S$.
%Let us now consider the class $[\hat Y]$
%of the divisor $\hat Y$ in ${\rm Pic}(P)=\CH^1(P)$.
There is an integer $m>0$ and a line bundle $\L$ on $S$ such that
$$
\O_{P}(\hat Y)\simeq {\rm pr}_1^*(\O_{{\bf P}^{d+1}_k}(m))\otimes_{\O_P}{\rm pr}_2^*(\L).
$$
If $d=0$,  take $H\subset {\bf P}^{1}_k$ to be the union of $m$ distinct points.
If $d>0$, take $H\subset {\bf P}^{d+1}_k$ to be a smooth irreducible hypersurface of degree $m$. Set $D=H\times_kS\subset P$ for the corresponding smooth ``horizontal" divisor 
in $P$. 
We can assume that $D$ intersects $\hat Y$ properly. Using Bertini's theorem, we can write $\L\simeq \O_S(T_1-T_2)$, where $T_1$ and $T_2$ are both smooth very ample divisors on $S$. 
Set $F_1={\rm pr}_2^{-1}(T_1)$, $F_2={\rm pr}_2^{-1}(T_2)$. 
Then we have
\begin{equation}\label{divf1f2}
\hat Y\sim D+F_1-F_2
\end{equation}
on $P$. %The divisor $D+F_1+F_2$ has strict normal crossings on $P$.
Let $U$ be an open subset of $P$ which contains the generic points 
of $D$, $F_1$, $F_2$, $\hat Y$, and is such that the divisor $U\cap (D+F_1+F_2+\hat Y)$ 
has strict normal crossings on $U$.
The work of Hironaka  on resolution of singularities now implies: 

\begin{thm}\label{Resol}
There is 
birational morphism $\beta: \ti P\to P={\bf P}^{d+1}\times_kS$ over $S$,
which is obtained as a succession of blow-ups along smooth centers
lying over $P-U$,
%(which at each step are transversal to the exceptional locus) 
such that the following is true: Let $T$ denote one of the divisors
$D$, $\hat Y$, $F_1$, $F_2$. Then for each such choice of $T$:

a) the strict transform $T'$ of  $T$ in $\ti P$ is smooth, 
%and the morphism $\beta_{|T'}: T'\to T$ is 
%a birational isomorphism,

b) the total transform $\beta^*(T)$ is a
divisor with
strict normal crossings and we can write
\begin{equation}\label{transform}
\beta^*(T)=T'+\sum_i n_i\cdot Z_i\ , 
\end{equation}
where $n_i\in \Z$ and the exceptional divisors $Z_i$   are birationally isomorphic over $S$ to projective bundles 
over smooth quasi-projective varieties of dimension $< \dim(X)$.
\end{thm}

\begin{Proof}
%As is also explained in [ML, Appendix A], 
This can be deduced from [H, Theorem I$^{N, n}_2$, pg. 170]
by taking $N=\dim(P)$, $n=N-1$, and $({\mathfrak R}_I^{N, N-1}, U)$
the resolution datum $((|D+F_1+F_2|; P; \hat Y), U)$.
(The last statement about the   divisors $Z_i$ follows
from the fact that, in this desingularization procedure,
the blow-up centers are always transverse to the exceptional 
locus.) 
 \endproof
\end{Proof}
\smallskip
 
By construction, $\ti P$ is a smooth variety that supports a projective morphism
$g: \ti P\to S$. By (\ref{divf1f2}) and Theorem \ref{Resol}, 
we now obtain
\begin{equation}\label{tra'}
Y'\sim D'+F_1'-F_2'+\sum_j m_j Z_j\ ,
\end{equation}
where  the  smooth divisors $Z_j$   are birationally isomorphic over $S$ to projective bundles over smooth quasi-projective varieties of dimension $<\dim(X)$.

Now let us use (i), (ii)$'$, (iii), (iv) 
of \S \ref{6b} (these hold under our induction hypothesis):
Apply the homomorphism $E^g: \Pic(\ti P)\to \CH^n(S)$ to the identity
(\ref{tra'}). Since the $Z_j$'s are smooth and  are birationally isomorphic over $S$ to projective bundles 
over smooth varieties of smaller dimension, (iii), (i), and (ii)$'$, give $E^g(Z_j)=E(Z_j\to S)=0$.
Using (iii) and (iv), we obtain $E^g( F_i')=Err(F_i'\to S)=0$, for $i=1$, $2$. On the other hand, 
 (ii)$'$  gives $E^g(D')=E(D'\to S)=E(D\to S)$. Thus, the identity (\ref{tra'}) implies
\begin{equation}\label{err2}
E^g(Y')=E(D\to S)\, .
\end{equation}
Now by its construction, $Y'$ is smooth and birationally isomorphic to $X$ over $S$. Hence, (iii) and (ii)$'$
imply $E^g(Y')=E(Y'\to S)=E(X\to S)$. Therefore, by (\ref{err2}),
to show $E(X\to S)=0$, it is enough to show that $E( D\to S )=0$.
Recall that $D=H\times_kS$. We  have
$
(pr_2)_*[\O_{H\times_kS}]=  \chi(H, \O_{H})\cdot [\O_S]\, 
$
in ${\rm K}_0(S)$. On the other hand, we find
$$
(pr_2)_*(\td_{n+d}(H\times_kS))=\frac{T_{d+n}}{T_d\cdot T_n}\cdot {\rm deg}(\td_d(H))\cdot \td_n(S)
$$
in $\CH^n(S)$. Therefore, since the Hirzebruch-Riemann-Roch for $H$ 
gives $ \chi(H, \O_{H})={\rm deg}(\td_d(H))/T_d$, we obtain $E(D\to S)=0$.
We conclude that integral Riemann-Roch theorem  for
$\F=\O_X$ holds for all $(f, n)$
with $d_f=d$. 

Now suppose
that $\F=\L$ is an invertible sheaf on $X$. 
We claim that
\begin{equation}
\frac{T_{d+n}}{T_n}\cdot\ctd_n(f_*[\L], S)=f_*(\ctd_{d+n}(\L, X ))\,.
\end{equation}
Set 
$\L=\O_X(-F)$, with $F$ a Weil divisor on $X$. Let us write
\begin{multline}\label{543}
\ctd_{n+d}(\L, X )=\ctd_{n+d}(\O_X, X )-\ctd_{n+d}([\O_X]-[\O_X(-F)], X )=\\ 
=\td_{n+d}( X )-\frac{T_{n+d}}{T_{n+d-1}}\cdot \td_{n+d}(F; X ) \ \ \ \ \ \ \ \ \ \ 
\end{multline}
(using (\ref{32}).) Apply Bertini's theorem (Proposition \ref{Bertini})
to the divisor $F$
on $X$. Using Proposition \ref{vample1} (c) for $m=n+d$ and (\ref{543}), 
we can express $\ctd_{n+d}(\L, X )$ in terms of  integral Todd classes of smooth varieties of dimension $\leq \dim(X)$. Correspondingly,
Proposition \ref{vample2} allows us to obtain a similar expression 
for the class $[\L]=[\O_X(-F)]$. Since, by the above, 
integral Riemann-Roch holds for the structure sheaf
when the relative dimension is $\leq d$, our claim for $\F=\L$  follows 
by 
comparing these two expressions.

Finally, suppose   that $\F$ is an arbitrary locally free coherent $\O_X$-sheaf on $X$.
By a result of Kleiman [K, Theorem 4.7 (b)], there is a birational  morphism $b: \ti X\to X$
which is obtained by successive  blow-ups along smooth centers, such that $ b^*\F$ 
has a filtration whose graded pieces are line bundles $\L_i$ on $\ti X$. In particular
\begin{equation*}
[b^*\F]=[\L_1]+\cdots +[\L_r]
\end{equation*}
in the Grothendieck group $\Kr_0(\ti X)$. Since both sides of the Riemann-Roch equation  are additive, our previous arguments give
that the integral Riemann-Roch identity holds for $b^*\F$ and $(f\cdot b, n)$. The argument in the proof of Corollary \ref{corfactor} 
now gives the integral Riemann-Roch identity
(\ref{RReqA}) for $\F$ and $(f, n)$. This and our induction, allows us to conclude the proof of Theorem \ref{main}.
\endproof
\medskip

%\begin{comment}
\begin{Remarknumb}\label{GenDeg}
{\rm 
Here we briefly sketch how a part of the argument 
can be recast using certain constructions 
in the theory of algebraic cobordism [LM]
and in particular the ``generalized degree formula". (We follow the notations of \S \ref{6b} and loc. cit.)

As a first step, one shows that the function $E: \SS_d\to \CH^n(S)$ given 
as in (\ref{errf}) gives a group homomorphism $E: \Omega^{-d}(S)\to \CH^n(S)$.
In fact, both terms of the difference (\ref{errf}) that defines 
the error function $E$ extend to group homomorphisms
from the cobordism group $\Omega^{-d}(S)$: For the first term, this follows
using the transformation between cobordism and $\rm K$-theory ([LM, 4.2]).
For the second, this can be seen using the transformation between 
cobordism and a suitable twisted Chow cohomology theory 
(see the argument in the construction
of the Conner-Floyd Chern 
classes [LM, 7.3.3], also [LM, 4.4.19]). 

Then, by the generalized degree formula [LM, Theorem 1.2.14],  we can find for $i=0,\ldots,  r$, classes $a_i\in \Omega^{-d_i}(k)$, and for $i=1,\ldots,  r$,  irreducible subvarieties $Z_i$ of $S$ together with morphisms $\ti Z_i \to S$ from smooth varieties $\ti Z_i$ with image $Z_i$  and $\ti Z_i\to Z_i$  birational, such that
\begin{equation}\label{gdf}
[f: X\to S]=a_0\cdot [{\rm id}: S\to S]+\sum_{i=1}^r a_i\cdot [\ti Z_i\to S]
\end{equation}
in $\Omega^{-d}(S)$. Here $d_0=d$ and $a_0$ is the degree of $[f]$ in the terminology of loc. cit.
In addition, by using embedded resolution for $Z_i\subset S$, 
the argument in the proof of [LM, Theorem 1.2.14] allows us to choose $\ti Z_i\to S$ that factor as a composition $\ti Z_i\hookrightarrow \ti S_i\to S$
where $\ti S_i$ is obtained by a successive blow-up of $S$ along smooth centers. We can also assume that for each $i$, the class $a_i$ is an integral linear  combination of the classes of smooth projective varieties $Y_{ij}\to \Spec(k)$.
Hence, (\ref{gdf}) implies that integral Riemann-Roch for $f$ and the structure sheaf will follow from
integral Riemann-Roch (for the structure sheaf) for $Y_{0j}\times S\to S$
and for the compositions 
$$
Y_{ij}\times \ti Z_i\xrightarrow{\, \rm pr\, } \ti Z_i\hookrightarrow \ti S_i\to S\ 
$$
with $i=1,\ldots , r$. Integral Riemann-Roch for a 
projection $Y\times Z\xrightarrow{\rm pr} Z$ now follows 
(as before) from Hirzebruch-Riemann-Roch for $Y\to \Spec(k)$ and the projection formula. On the other hand,  integral Riemann-Roch for 
$\ti Z_i\hookrightarrow \ti S_i$ and $\ti S_i\to S  $
can be shown using Theorem \ref{JoRR},  and Theorem \ref{blowRR} and Proposition \ref{compose} respectively. Integral Riemann-Roch for $(f, n)$ and the structure sheaf now follows from the above and Proposition 
\ref{compose}.

}
\end{Remarknumb} 
%\end{comment}
 \bigskip
 
 % \vfill\eject

\end{document}